\documentclass{article}

\usepackage{microtype}
\usepackage{graphicx}
\usepackage{subcaption}
\usepackage{booktabs} 

\usepackage{amsmath}
\usepackage{amsthm}
\usepackage{amssymb}
\usepackage{physics}
\usepackage{float}

\newtheorem{theorem}{Theorem}
\newtheorem{lemma}{Lemma}

\newtheorem{definition}{Definition}
\newtheorem{assumption}{Assumption}

\newcommand{\R}{\mathbb{R}}

\newcommand{\N}{\mathbb{N}}
\newcommand{\deq}{\mathrel{\mathop:}=}

\usepackage{hyperref}

\usepackage[accepted]{icml2020}

\icmltitlerunning{Randomised Subspace Gauss-Newton (R-SGN)}

\begin{document}

\twocolumn[
\icmltitle{A Randomised Subspace Gauss-Newton Method for Nonlinear Least-Squares}



\icmlsetsymbol{equal}{*}

\begin{icmlauthorlist}
\icmlauthor{Coralia Cartis}{Oxford,ATI,equal}
\icmlauthor{Jaroslav Fowkes}{Oxford,equal}
\icmlauthor{Zhen Shao}{Oxford,equal}
\end{icmlauthorlist}

\icmlaffiliation{Oxford}{Mathematical Institute, University of Oxford, Woodstock Road, Oxford OX2 6GG, UK.}
\icmlaffiliation{ATI}{The Alan Turing Institute for Data Science, The British Library, London NW1 2DB, UK}

\icmlcorrespondingauthor{Jaroslav Fowkes}{jaroslav.fowkes@stfc.ac.uk}

\icmlkeywords{optimization,least-squares}

\vskip 0.3in
]
\textbf{}


\printAffiliationsAndNotice{*The order of authors is alphabetical.} 

\begin{abstract}
We propose a Randomised Subspace Gauss-Newton (R-SGN) algorithm for solving nonlinear least-squares optimization problems,
that uses a sketched Jacobian of the residual in the variable domain and solves a reduced linear least-squares on each iteration.
A sublinear global rate of convergence result is presented for a trust-region variant of R-SGN, with high probability, which matches deterministic counterpart results in the order of the accuracy tolerance.
Promising preliminary numerical results are presented for R-SGN on logistic regression and on nonlinear regression problems from the CUTEst collection.
\end{abstract}

\section{Introduction}
We aim to solve the nonlinear least-squares problem
\begin{align}
\min_{x \in \R^d} f(x) = \frac{1}{2}\sum_{i=1}^n \norm{ r_i(x) }_2^2=\frac{1}{2}\|r(x)\|_2^2  \label{NLS}
\end{align}
where $r = (r_1, \dots, r_n): \R^d \to \R^n$ is a smooth vector of nonlinear (possibly nonconvex) residual functions. We define the Jacobian (matrix of first order derivatives) as
\begin{align*}
J(x) = \left(\pdv{r_i(x)}{x_j}\right)_{ij} \in \R^{n \times d}
\end{align*}
and can then compactly write the gradient as $\grad f(x) = J(x)^Tr(x)$.
The classical Gauss-Newton (GN) algorithm simply applies Newton's method to minimising $f$ with only the first-order $J(x)^TJ(x)$ term in the Hessian, dropping the second-order terms involving the Hessians of the residuals $r_i$. 
Thus, at every iterate $x_k$, Gauss-Newton approximately minimises the following convex quadratic local model
\begin{align*}
f(x_k) + \langle J(x_k)^Tr(x_k), \hat{s} \rangle + \frac{1}{2} \langle \hat{s}, J(x_k)^T J(x_k) \hat{s} \rangle 
\end{align*}
over $\hat{s} \in \R^d$. In our approach, which we call Randomised Subspace Gauss-Newton (R-SGN), we reduce the dimensionality of this model by minimising in an $l$-dimensional randomised subspace $\mathcal{L} \subset \R^d$, with $l\ll d$, by approximately minimising the following reduced model 
\begin{align}
f(x_k) + \langle J_{\mathcal{S}}(x_k)^Tr(x_k), s \rangle + \frac{1}{2} \langle s, J_{\mathcal{S}}(x_k)^T J_{\mathcal{S}}(x_k) s \rangle 
\label{m_k_definition_TR}
\end{align}
over $s \in \R^l$, where  $J_{\mathcal{S}}(x_k) \in \R^{n\times l}$ denotes the reduced Jacobian. Compared to the classical Gauss-Newton model, this reduced model also offers the computational advantage that it only needs to evaluate $l$ Jacobian actions, giving $J_{\mathcal{S}}(x_k)$, instead of the full Jacobian matrix $J(x_k)$. 

In its simplest form, $J_{\mathcal{S}}$ can be thought of as a random subselection of columns of the full Jacobian $J$, which leads to variants of our framework that are Block-Coordinate Gauss-Newton (BC-GN) methods. In this case, for example, if the Jacobian were being calculated by finite-differences of the residual $r$, only a small number of evaluations of $r$ along coordinate directions would be needed;  such a BC-GN variant has already been used for parameter estimation in climate modelling \cite{tett2017calibrating}. 

\section{Related Work and Motivation}
Block-coordinate descent methods have been, by now, extensively studied, especially for minimizing convex functions \cite{richtarik2014iteration}, but not only;
the randomised nonmonotone block proximal-gradient method \cite{lu2017randomized} minimises the sum of a smooth (possibly nonconvex) function and a block-separable (possibly nonconvex nonsmooth) function, with global rates of convergence being provided for when the expected value of the (true) gradient is sufficiently small. In a similar vein, \citet{facchinei2015parallel} propose a general decomposition framework for the parallel optimization of the sum of a smooth (possibly nonconvex) function and a block-separable nonsmooth convex function. Both of these proposals are gradient descent frameworks, while for improved and robust practical performance,
(some) second-order information of the objective function needs to be employed in the algorithm.  As second-order methods may be too computationally expensive for large-scale applications, subspace variants have been devised; extensive literature exists on deterministic choices such as Newton-CG (conjugate gradient) and Krylov methods, both for calculating approximate Newton directions for general optimization and approximate Gauss-Newton directions for nonlinear least-squares \citep[see for example,][]{nocedal2006numerical, gratton2007approximate}. However, though these methods are really powerful and very much state-of-the-art for inexact second-order methods, they still require (full) Hessian-vector products, which may be too expensive for some applications \cite{tett2017calibrating}.
Also, we are interested in using randomisation techniques for choosing the subspace of minimization, in an attempt to exploit the benefits of Johnson-Lindenstrauss (JL) Lemma-like results, that essentially reduce the dimension of the optimization problem without loss of information.
Such sketching techniques have already been proposed, especially for convex optimization, as we describe next.
 The sketched Newton algorithm \cite{pilanci2017newton} requires a sketching matrix that is proportional to the rank of the Hessian, which may be too computationally expensive. By contrast, sketched online Newton \cite{luo2016efficient}  uses streaming sketches to scale up a second-order method, comparable to Gauss–Newton, for solving online learning problems.
The randomised subspace Newton \cite{gower2019rsn} efficiently sketches the full Newton direction for a family of generalised linear models, such as logistic regression.
The stochastic dual Newton ascent algorithm in  \citet{qu2016sdna} requires a positive definite upper bound $M$ on the Hessian and proceeds by selecting random principal submatrices of $M$ that are then used to form and solve an approximate Newton system. 
The randomized block cubic Newton method  in \citet{doikov2018randomized} combines the ideas of randomized coordinate descent with cubic regularization and requires the optimization problem to be block separable.

Here, we propose a randomized subspace Gauss-Newton method for nonlinear least-squares problems, that, at each iteration,  only needs a sketch of the Jacobian matrix in the variable domain, which it then uses to solve a reduced linear least-squares problem for the step calculation. To ensure global convergence of the method, from any starting point, we include a trust-region technique (and also have results for a quadratic regularization variant); unlike most prior work, we focus on the general nonconvex case of problem \eqref{NLS}.

\section{The R-SGN Algorithm}
As mentioned above, at each iterate $x_k$, $k\geq 0$, R-SGN uses Jacobian actions $J_{\mathcal{S}}(x_k)$ in an $l$-dimensional randomised subspace $\mathcal{L} \subset \R^d$, $l \ll d$. To obtain these actions we draw a random sketching matrix $S_k \in \R^{l \times d}$ from a chosen random matrix distribution $\mathcal{S}$ at each iteration, and set 
\[
J_{\mathcal{S}}(x_k) \deq J(x_k)S_k^T
\]
in the local model \eqref{m_k_definition_TR}.   
This model is then approximately minimized over a trust region ball, $\norm{s}_2 \leq \Delta_k$, so that a typical 
 Cauchy decrease condition \citep[see e.g.][]{conn2000trust} is satisfied,
 \begin{align}
m_k(0) - m_k(s_k) \geq c_1\|g_k\|_2 \min \left( \Delta_k, \frac{\|g_k\|_2}{\|B_k\|_2} \right) \label{m_k_termination_condition},
\end{align}
for some constant $c_1 >0$; where $g_k\deq J_{\mathcal{S}}(x_k)^Tr(x_k)$
and $B_k\deq J_{\mathcal{S}}(x_k)^TJ_{\mathcal{S}}(x_k)$.
The step is then either accepted or rejected depending on the ratio $\rho_k$ between the actual objective decrease and the decrease predicted by the model, 
\begin{align}
\rho_k \deq \frac{f(x_k)-f(x_k + S_k^T s_{k})}{f(x_k) - m_k(s_{k})}. \label{eq:rho}
\end{align}
A full description of the R-SGN with trust region algorithm is provided in \autoref{alg:rsgn} below.

\begin{algorithm}[tbh]
   \caption{R-SGN with Trust Region}
   \label{alg:rsgn}
\begin{algorithmic}
   \STATE {\bfseries Input:} $x_0 \in \R^{d}$
   \STATE {\bfseries Parameters: } $l<d$, $\eta \in (0,1)$, $0 < \gamma_1 <1 < \gamma_2$, such that $\gamma_2 = \gamma_1^{-c}$ for some $c\in \N$, $\Delta_0 \in \R$.
   \FOR{$k=0$ {\bfseries to} $N$}
   \STATE Randomly draw a sketching matrix $S_k \in \R^{l\times d}$.
   \STATE{ Compute step $s_k$ to approximately minimise the model \eqref{m_k_definition_TR} satisfying \eqref{m_k_termination_condition}. }
		
    \STATE{ Compute decrease $\rho_k$ given by \eqref{eq:rho}.} 
    \IF{$\rho_k \geq \eta$}
    \STATE{$x_{k+1} = x_k + S_k^T s_{k}$ }
    \STATE{$\Delta_{k+1} = \gamma_2 \Delta_k$}
    \ELSE
    \STATE{$\Delta_{k+1} = \gamma_1 \Delta_k$}
    \ENDIF
\ENDFOR
\end{algorithmic}
\end{algorithm}

\section{A global rate of convergence for R-SGN}






We assume the following properties of the residual $r(x)$, the Jacobian $J(x)$ and the random sketching matrices $S_k$.




\begin{assumption}\label{AA1}
The residual $r$ is continuously differentiable, and $r$ and its Jacobian $J$ are Lipschitz continuous on $\R^d$. 
\end{assumption}

\begin{assumption}\label{AA2}
Let $\epsilon_S, \delta_S \in (0,1)$. At each iterate $x_k$, with probability at least $1-\delta_S$, we have that 
\begin{align}
\|S_k \grad f(x_k)\|_2^2 \geq (1-\epsilon_S) \|\grad f(x_k)\|_2^2, \label{eqn::partialJL}
\end{align}
and \begin{align}
\|S_k\|_2 \leq S_{max}. \label{eqn::S_k_bound}
\end{align}
\end{assumption}
\begin{assumption}\label{AA3}
The probability in \autoref{AA2} satisfies
\begin{align}
1- \delta_S > c_2,
\end{align}
where $c_2 = \frac{c+2}{2c+2}$, with $c$ being a user chosen parameter defined in \autoref{alg:rsgn}. If $c=1$, i.e.\ we have $\gamma_2 = \gamma_1^{-1}$ in \autoref{alg:rsgn}, then $c_2 = 3/4$.
\end{assumption}
The following result shows that the R-SGN algorithm produces an iterate with arbitrarily small $\|\grad f(x_k)\|_2$, with exponentially small failure probability, in a quantifiable number of iterations. 

	\begin{theorem}\label{thm::TR::complexity}
	Let Assumptions \ref{AA1}, \ref{AA2}, \ref{AA3} hold. Let $\epsilon>0$, and  $\delta_1 \in (0,1)$ 
	such that $(1-\delta_S)(1-\delta_1) > c_2$. Then the R-SGN algorithm takes at most 
	\begin{align}
 N\leq \frac{1}{(1-\delta_S)(1-\delta_1)-c_2}\mathcal{O} \left( \frac{ f(x_0) - f^* }{(1-\epsilon_S) \epsilon^2} \right) \label{eq::thm::QR::complexity}
	\end{align}	
	iterations and evaluations of the residual and sketched Jacobian such that $\min_{k \leq N} \|\grad f(x_k)\|_2 \leq \epsilon $, with probability at least $1-e^{-\frac{\delta_1^2}{2} (1-\delta_S) N}$.	
	\end{theorem}
We note that as a function of the tolerance $\epsilon$, this bound matches deterministic complexity bounds for first-order and Gauss-Newton  methods 
\cite{cartis2012much}, despite having only partial Jacobian information available at each iteration.

We also have similar global rates of convergence results for R-SGN with quadratic regularization techniques instead of trust region.
There is no convexity or special structure requirement on $f$ or $r$ in our results. 
 Our proof employs techniques from convergence/complexity analysis of optimisation methods based on probabilistic models in \citet{Cartis:2017fa,Gratton:2017kz}. However, our result only assumes a weaker condition, namely, that the 'model gradient' $S_k \grad f(x_k)$ has similar norm to that of the full gradient $\grad f(x_k)$, while previous results assume the model gradient itself approximates the true gradient. 
\subsection{Discussion on generating the random matrix $S_k$}

We give a few suitable choices of the random matrix ensemble $\cal{S}$ satisfying \autoref{AA2}. 
\begin{definition}[Non-uniformity dependent oblivious JL embedding]\label{def::1}
Let $\epsilon_S, \delta_S \in (0,1)$, $\nu_S \in (0,1]$, $d \in \N$. A distribution on $S \in \R^{l \times d}$ with $l = l(d, \nu_S, \epsilon_S,\delta_S)$ is a $(\nu_S, \epsilon_S,\delta_S)$-oblivious JL embedding if, given fixed $y \in \R^{d}$ with $\|y\|_{\infty}/ \|y\|_2 \leq \nu_S$, we have that with probability at least $1-\delta_S$, a matrix $S$ drawn from the distribution satisfies
\begin{align}
(1-\epsilon_S)\|y\|_2^2 \leq \|Sy\|_2^2 \leq (1+\epsilon_S) \|y\|_2^2.
\end{align}
\end{definition}
Thus  we can satisfy \eqref{eqn::partialJL} with probability at least $1-\delta_S$ at each R-SGN iteration if our random matrix distribution $\cal{S}$ is a $(\nu_S,\epsilon_S, \delta_S)$-oblivious JL embedding, and we have that $\|\grad f(x_k)\|_{\infty} / \|\grad f(x_k)\|_2 \leq \nu_S$ for all iterates $x_k$ by taking $y = \grad f(x_k)$ in \autoref{def::1}.\footnote{$\|\grad f(x_k)\|_{\infty} / \|\grad f(x_k)\|_2 \leq \nu_S$ is not required if $\nu_S=1$, as then the inequality always holds.}

The well studied oblivious JL-embedding \citep[see][]{Woodruff:2014aa}, is a special case of \autoref{def::1} with $\nu_{S}=1$. In particular, we have the following result.
\begin{lemma}[\citealp{Dasgupta:2002aa}]
We say $S \in \R^{ l \times d }$ is a Gaussian matrix if $S_{ij}$ are identically distributed as $N (0, {l}^{-1/2})$. Let $\epsilon_S , \delta_S \in (0,1)$, then the distribution of Gaussian matrices $S \in \R^{l \times d}$ with $l = \mathcal{O}(|\log\delta_S|/\epsilon_S^2)$ is an $(1, \epsilon_S, \delta_S )$-oblivious JL embedding. 
\end{lemma}
The following hashing/sparse embedding matrices are also $(1,\epsilon_S,\delta_S)$-oblivious JL embeddings.
\begin{lemma}[\citealp{Kane:2012aa,Cohen:2018fi}]
We define $S \in \R^{l \times d}$  to be a s-hashing matrix if, independently for each $j \in [d]$, we sample without replacement $i_1, i_2, \dots, i_s \in [l]$ uniformly at random and let $S_{i_k j} = \pm 1/\sqrt{s}$, $k = 1, 2, \dots, s$. Let $\epsilon_S , \delta_S \in (0,1)$, then the distribution of s-hashing matrices $S \in \R^{l \times d}$ with $s = \mathcal{O}(|\log\delta_S|/\epsilon_S) $ and $l = 
\mathcal{O}\left(|\log \delta_S|/ \epsilon_S^2  \right) $ is an $(1, \epsilon_S, \delta_S )$-oblivious JL embedding.
\end{lemma}
The distribution of sampling matrices is not usually an oblivious JL embedding. However it can be shown to be a non-uniformity dependent oblivious JL embedding. 
\begin{lemma}\label{lem:sampling}
We define $S \in \R^{l \times d}$ to be a sampling matrix if, independently for each $i \in [l]$, we sample $j \in [d]$ uniformly at random and let $S_{ij}=(d/l)^{-1/2}$. Let $\epsilon_S , \delta_S, \nu_S \in (0,1)$, then the distribution of sampling matrices $S \in \R^{l \times d}$ with $l = \mathcal{O}(d\nu_S ^2|\log\delta_S|/\epsilon_S^2) $ is an $( \epsilon_S, \delta_S, \nu_S )$-oblivious JL embedding.
\end{lemma}
Note that R-SGN with a sampling choice for $S_k$  leads to  block-coordinate variants (BC-GN) of R-SGN, with fixed block size $l$.

The above properties of Gaussian and hashing matrices imply that their embedding dimension $l$ is independent of the variable dimension $d$, while for sampling matrices $l \propto d\nu_S^2$; this has direct implications on the size of the sketch in R-SGN; in particular, it does not have to grow asymptotically. However, when $\| \grad f(x_k)\|_{\infty} \leq C d^{-1/2}\| \grad f(x_k)\|_2$ for some constant $C$ (thus
the maximum gradient component is no larger than a constant multiple of the average gradient component), then sampling matrices also have an embedding dimension $l$ that is independent of the variable dimension $d$. Thus the more uniformly distributed the magnitude of the gradient entries, the better sampling works in R-SGN (namely, in BC-GN), ensuring a.s.\ global convergence. 

It can be easily verified that \eqref{eqn::S_k_bound} is satisfied for sampling and hashing matrices. For Gaussian matrices, one may use  the non-asymptotic bound established for the maximum singular value of Gaussian matrices in \citet{rudelson2010non} to show $\|S_k\|_2$ is bounded with high probability. \autoref{AA2} can then be satisfied with a union bound.

\section{Preliminary Numerical Results}
For simplicity, in our performance evaluation of R-SGN, we use sampling sketching matrices as defined in \autoref{lem:sampling},
which leads to a BC-GN variant of R-SGN.
We first consider logistic regression, written in the form \eqref{NLS}, by letting
$r_i(x) = \ln(1 + \exp(-y_i a_i^T x))$,
where $a_i \in \R^d$ are the observations and $y_i \in \{-1,1\}$ are the class labels; we also include a quadratic regularization term $\lambda\|x\|_2^2$ by treating it as
an additional residual.

\begin{figure}[!t]
\begin{subfigure}{\columnwidth}
\includegraphics[width=\columnwidth]{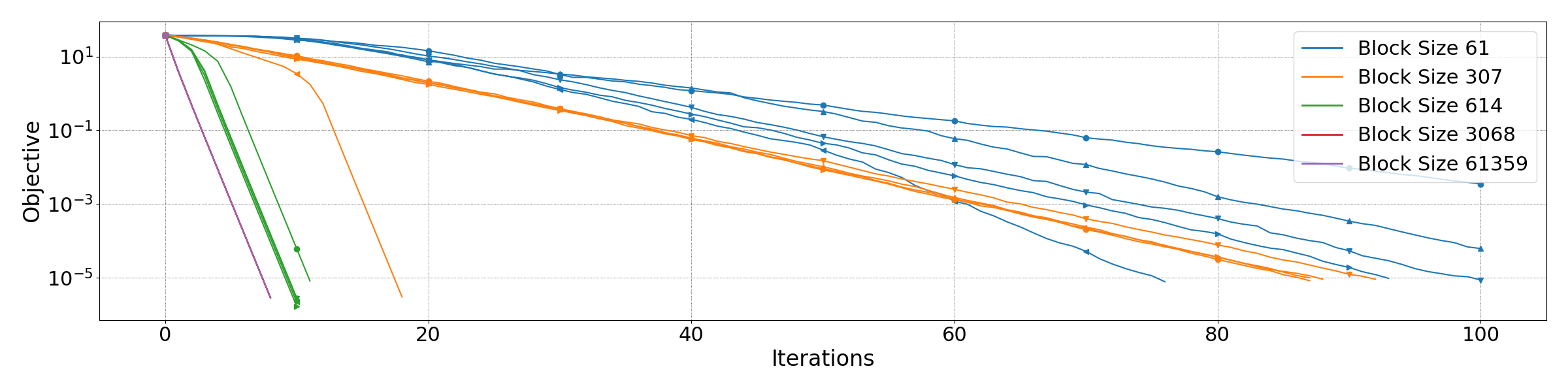}
\end{subfigure}
\begin{subfigure}{\columnwidth}
\includegraphics[width=\columnwidth]{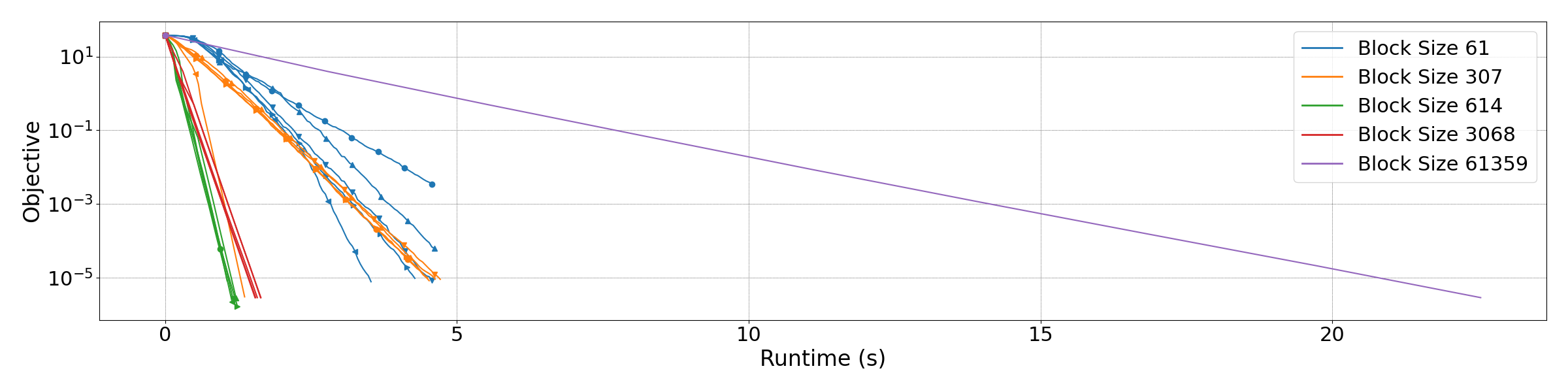}
\end{subfigure}
\caption{R-SGN on the \textsc{chemotherapy} dataset}\label{fig:chemo}
\end{figure}

\begin{figure}[!t]
\begin{subfigure}{\columnwidth}
\includegraphics[width=\columnwidth]{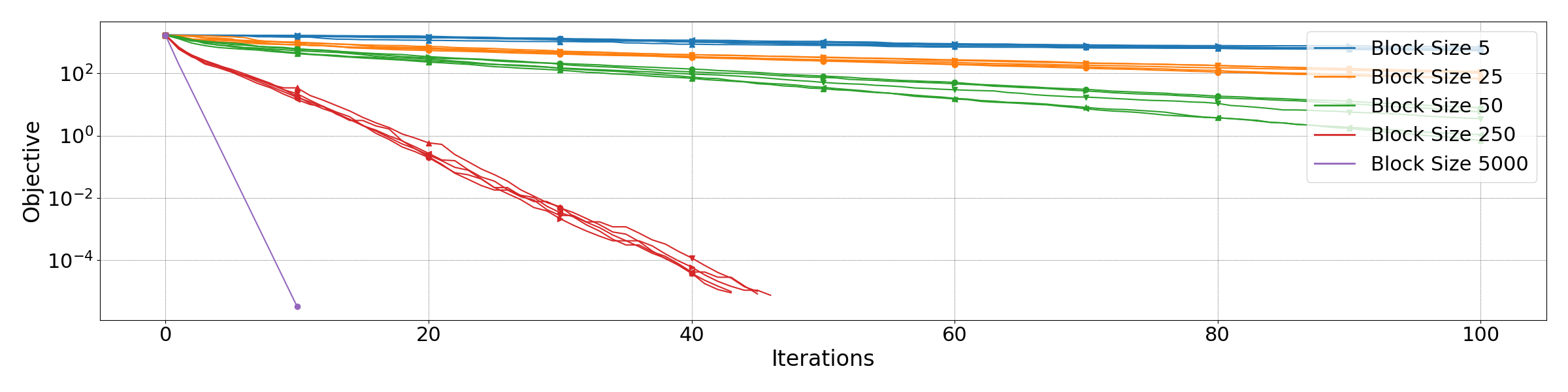}
\end{subfigure}
\begin{subfigure}{\columnwidth}
\includegraphics[width=\columnwidth]{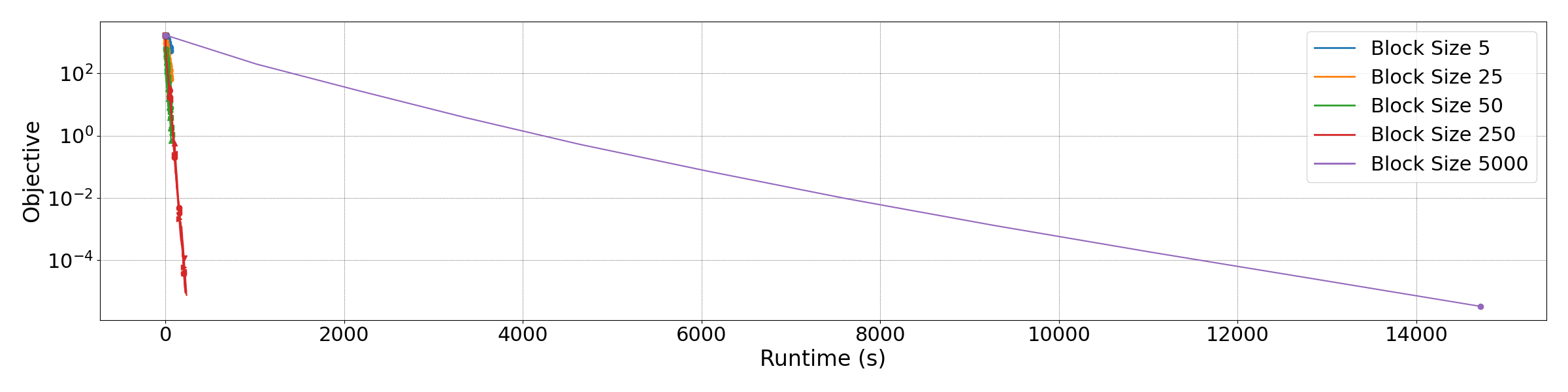}
\end{subfigure}
\caption{R-SGN on the \textsc{gisette} dataset}\label{fig:gisette}
\end{figure}

We first test on the \textsc{chemotherapy} and \textsc{gisette} datasets from OpenML \cite{OpenML2013} for $100$ iterations with $\lambda=10^{-10}$, using block sizes of 0.1\%, 0.5\%, 1\%, 5\% and 100\% of the original for the $61,359$ dimensional \textsc{chemotherapy} dataset and the $5,000$ dimensional \textsc{gisette} dataset; in a similar testing setup to \citet{gower2019rsn}. As R-SGN is randomised, we perform five runs of the algorithm for each block size starting at $x_0=0$. We terminate once the objective $f(x_k)$ goes below $10^{-5}$ and plot $f(x_k)$ against iterations and runtime in each Figure. On the \textsc{chemotherapy} dataset, we see from \autoref{fig:chemo} that we are able to get comparable performance to full GN ($d=61,359$ in purple) using only $1\%$ of the original block size ($l=614$ in green) at $1/20$th of the runtime.  For the \textsc{gisette} dataset, we see from \autoref{fig:gisette} that similarly, we are able to get good performance compared to GN ($d=5,000$ in purple) using $5\%$ of the original block size ($l=250$ in red) at $1/60$th of the runtime. 

\begin{figure}[!t]
\begin{subfigure}{\columnwidth}
\includegraphics[width=\columnwidth]{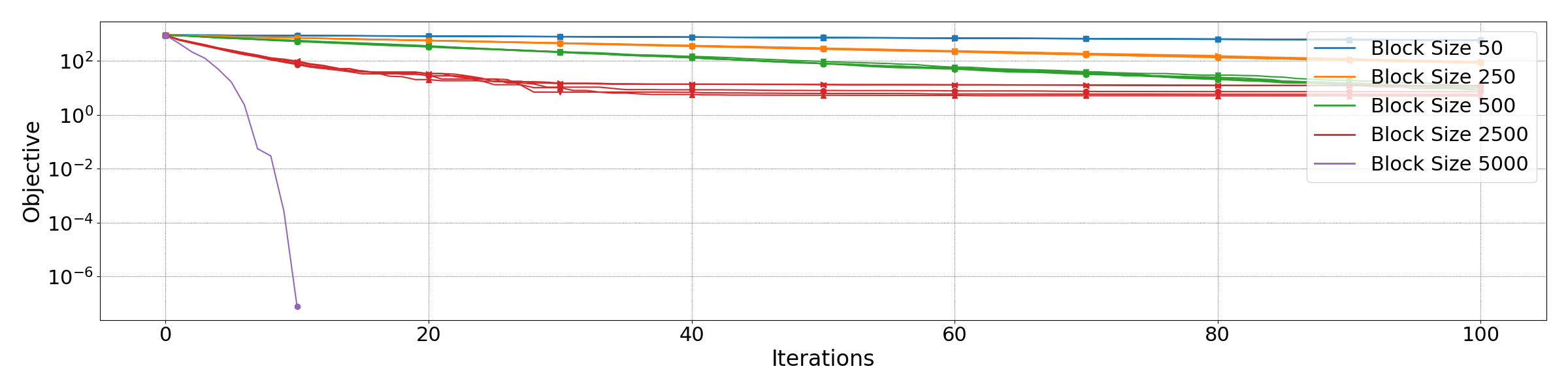}
\end{subfigure}
\begin{subfigure}{\columnwidth}
\includegraphics[width=\columnwidth]{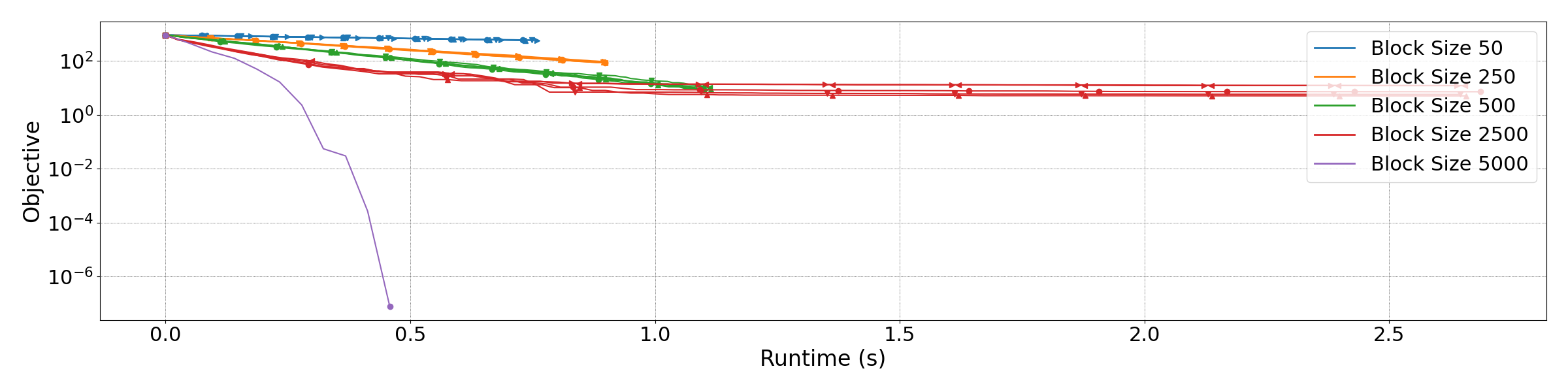}
\end{subfigure}
\caption{R-SGN on the \textsc{artif} dataset}\label{fig:artif}
\end{figure}

\begin{figure}[!t]
\begin{subfigure}{\columnwidth}
\includegraphics[width=\columnwidth]{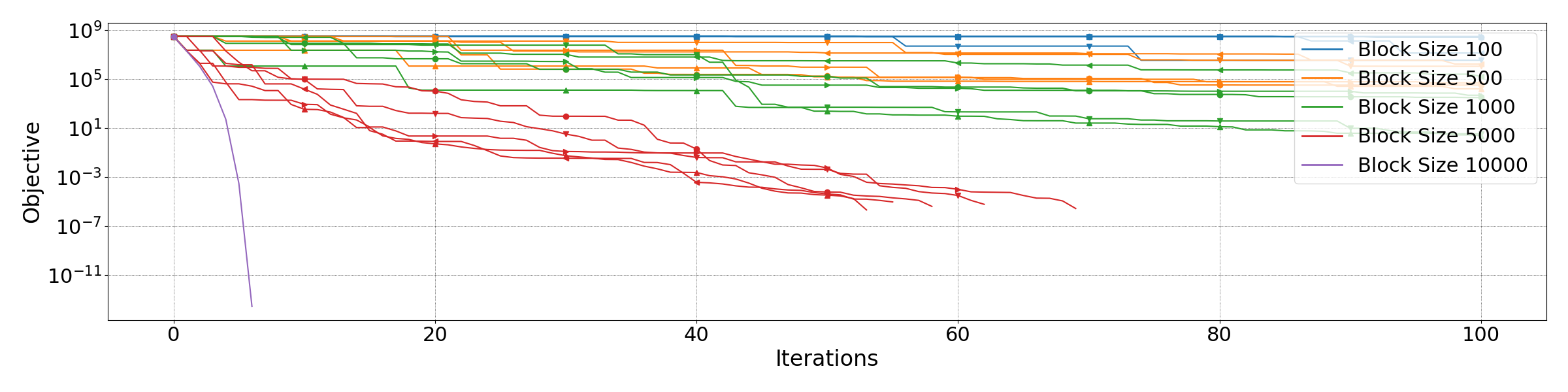}
\end{subfigure}
\begin{subfigure}{\columnwidth}
\includegraphics[width=\columnwidth]{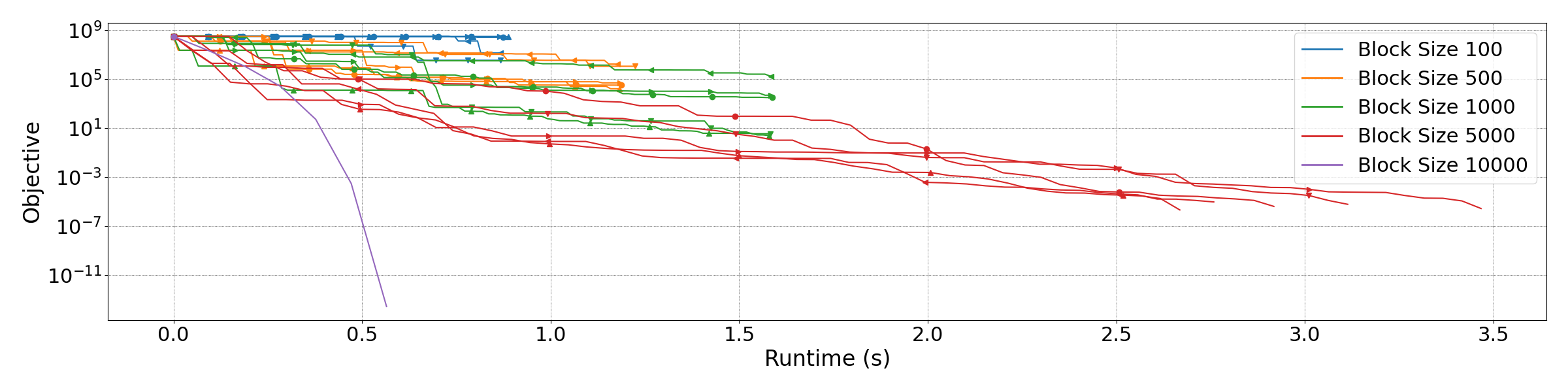}
\end{subfigure}
\caption{R-SGN on the \textsc{oscigrne} dataset}\label{fig:oscigrne}
\end{figure}

We next test on the \textsc{artif} and \textsc{oscigrne} nonlinear least-squares problems from the CUTEst collection \cite{gould2015cutest} with dimension $5,000$ and $10,000$ respectively in both $n$ and $d$.  As before we run R-SGN five times for each problem for $100$ iterations from the given starting points for each problem, with block sizes of 1\%, 5\%, 10\%, 50\% and 100\% of the original. As we can see from the results in \autoref{fig:artif} and \autoref{fig:oscigrne}, R-SGN does not perform as well as on the logistic regression problems, due possibly to the increased difficulty/nonlinearity of these problems. However, it is still able to reach low accuracy solutions reasonably well, achieving a reduction of several orders of magnitude in the objective after just one second for both problems at $50\%$ block size (plotted in red). 
The CUTEst examples show that more sophisticated approaches may be needed to tackle general nonlinear least-squares problems; for example instead of keeping the block size fixed at $l$, we may consider increasing $l$ over iterations, a topic of future investigation.

\vspace*{-0.25cm}
\section{Conclusions}
We have presented an optimization algorithm for solving nonlinear least-squares problems that uses a sketched Gauss-Newton approximation to the objective's second derivatives, in a randomised subspace of the full space. A global rate of convergence, with high probability, is presented. We have demonstrated promising numerical results on logistic regression problems and nonlinear CUTEst test problems.  

\bibliography{example_paper}
\bibliographystyle{icml2020}

\end{document}